\documentclass[11pt]{article}
\usepackage{enumitem}
\usepackage{amsmath,amsthm,amssymb,amsfonts}
\usepackage{amsmath}
\usepackage{multirow}
\usepackage{graphicx}
\usepackage[colorlinks=true,linkcolor=black]{hyperref}
\usepackage{hyperref}
\usepackage{cancel}
\usepackage[toc,page]{appendix}
\usepackage{fullpage}
\usepackage[margin=0.75in]{geometry}
\usepackage[utf8]{inputenc}

\theoremstyle{plain}
\newtheorem{theorem}{Theorem}[section]

\theoremstyle{definition}

\newtheorem{remark}[theorem]{Remark}
\newtheorem*{theorem*}{Theorem}

\numberwithin{equation}{section}

\def\R{\mathbb{R}}

\def\S{\mathbb{S}^{n-1}}

\def\t{\nabla_{T}}

\def\L2{L^2(\S)}

\def\s{\mathbb{S}^2}
\def\ds{\dashint_{\mathbb{S}^{n-1}}}
\def\d2{\dashint_{\mathbb{S}^2}}

\def\Xint#1{\mathchoice
{\XXint\displaystyle\textstyle{#1}}%
{\XXint\textstyle\scriptstyle{#1}}%
{\XXint\scriptstyle\scriptscriptstyle{#1}}%
{\XXint\scriptscriptstyle\scriptscriptstyle{#1}}%
\!\int}
\def\XXint#1#2#3{{\setbox0=\hbox{$#1{#2#3}{\int}$ }
\vcenter{\hbox{$#2#3$ }}\kern-.6\wd0}}

\def\dashint{\Xint-}

\def\ds{\dashint_{\S}}

\begin{document}

\centerline{{\Large {\bf A note on a rigidity estimate for degree $\pm 1$ conformal maps on $\s$}}}

\vspace*{0.5  cm}

\centerline{{\large {\bf Jonas Hirsch$^1$,\ \  Konstantinos Zemas$^2$}}}
 
\vspace*{0.65  cm}

\begin{abstract}
\footnotesize In this note we present a short alternative proof of an estimate obtained by A.B.-Mantel, C.B. Muratov and T.M. Simon in \cite{bernand2019quantitative} regarding the rigidity of degree {$\pm 1$} conformal maps of $\s$, i.e. its Möbius transformations.\\[-10pt]
\end{abstract}

\vspace{0.3 cm}


\section{Introduction}\label{sec:1}
As it is well known in the theory of harmonic mappings, a map from $\mathbb{S}^2$ into a Riemannean manifold is harmonic, i.e. a critical point of the Dirichlet energy functional iff it is generalized conformal. In particular, according to the famous theorem by J. Liouville, the class of orientation-preserving (respectively orientation-reversing) degree 1 (respectively -1) harmonic maps from $\s$ onto itself is precisely the group of its Möbius transformations. 

In this note we would like to show how a proof of Liouville's theorem that is included in Appendix A of \cite{zemas2020} can be perturbed quantitatively to give an alternative, shorter proof of a result obtained by A.B.-Mantel, C.B. Muratov and T.M. Simon in \cite{bernand2019quantitative}. In particular, Theorem 2.4 therein  can be thought of as a sharp quantitative version of Liouville's theorem regarding degree {$\pm 1$} conformal maps from $\s$ to itself.\\[-12pt]

We discuss here the case of maps of degree $1$, the case of maps of degree $-1$ being completely analogous, or it can simply be derived by the previous one after composing with the orientation-reversing orthogonal map $x:=(x_1,x_2,x_3)\in \s\mapsto \tilde{x}:=(x_1,x_2,-x_3)\in \s$. Let us denote by 
\begin{equation}\label{admissible_class_s_2}
\mathcal{A}_{\s}:=\left\{u\in W^{1,2}(\s;\s):\  \mathrm{deg}u:=\dashint_{\s}\big\langle u,\partial_{\tau_1}u\wedge\partial_{\tau_2}u\big\rangle=1\right\}
\end{equation}
and refer the reader to \cite{brezis1995degree} for more details on the notion of degree for Sobolev maps among closed Riemannean manifolds of the same dimension. Here we have used the standard ``Hodge-dual'' identification
\begin{equation*}
\big\langle u,\partial_{\tau_1}u\wedge\partial_{\tau_2}u\big\rangle:=\big\langle u,\ast(\partial_{\tau_1}u\wedge\partial_{\tau_2}u)\big\rangle,
\end{equation*}
where $\{\tau_1, \tau_2\}$ is a local orthonormal frame for $\s$, indicated by the unit normal vector. Let us also denote by $Conf(\s)$ the group of Möbius transformations of $\s$ and by $Conf_+(\s)$ the subgroup of those that are orientation-preserving.
As it is well known, the excess energy
\begin{equation}\label{deficit_on_S2}
D(u):=\frac{1}{2}\d2|\t u|^2-1
\end{equation}
considered among maps in $\mathcal{A}_{\s}$, is nonnegative, invariant under precompositions with elements of $Conf_+(\s)$, vanishes precisely when $u\in Conf_+(\s)$ and therefore provides an appropriate notion of conformal deficit for maps in $\mathcal{A}_{\s}$.\\[-12pt] 

Here we have also adopted the convention that $\t u$ denotes the extrinsic gradient of $u\in W^{1,2}(\s;\s)$, i.e. the gradient of $u$ when this is considered as a $W^{1,2}$-map from $\s$ to $\R^3$ with $|u|\equiv 1$ $\mathcal{H}^2$-a.e. on $\s$. With this notation, one has the following result (see \cite{bernand2019quantitative}, Theorem 2.4).
\begin{theorem}\label{main_thm}
There exists a universal constant $c>0$ so that for every $u\in \mathcal{A}_{\s}$ there exists $\phi \in Conf_+(\s)$ such that
\begin{equation}\label{conf_s2_s2_quantitative}
\d2 \left|\t u-\t \phi\right|^2\leq c D(u).\\[10pt]
\end{equation}
\end{theorem}

Let us also mention that another proof of the previous theorem was given recently by P. Topping in \cite{topping2020}, using the harmonic map heat flow on $\s$. We comment briefly on the latter approach at the end of the note.

\section{Proof of \hyperref[main_thm]{Theorem 1.1.}}
{\bf{Step 1}}. By an easy topological argument that is explained also in the proof of Liouville's theorem in \cite{zemas2020}, \textit{given $u\in \mathcal{A}_{\s}$ one can always find $\psi\in Conf_+(\s)$ so that $\d2 u\circ\psi=0$}.
Let us revise the argument for this.\\[-12pt]

Assume first that $u\in \mathcal{A}_{\s}\cap C^{\infty}$, in particular $u$ is surjective. If $b_u:=\d2 u=0$ there is nothing to prove. If $b_u\neq 0$, one can show that there always exist $\xi_0 \in \s$ and $\lambda_0>0$ so that
\begin{equation*}\label{Möbius_transformation_fixing_the_center}
\d2 u\circ\phi_{\xi_0,\lambda_0}=0,
\end{equation*}
where $\phi_{\xi_0,\lambda_0}:=\sigma_{\xi_0}^{-1}\circ i_{\lambda_0}\circ\sigma_{\xi_0}$, with $\sigma_{\xi_0}$ being the stereographic projection of $\s$ onto $T_{\xi_0}\s\cup\{\infty\}$ and \linebreak 
$i_{\lambda_0}:T_{\xi_0}\s\mapsto T_{\xi_0}\s$ being the dilation in $T_{\xi_0}\s$ by factor $\lambda_0>0$.\\[-12pt] 

Indeed, consider the map $F:\s \times [0,1]\mapsto \overline {B_1}$ defined as 
\begin{equation*}\label{homotopy_for_fixing_the_center}
F(\xi,\lambda):=\d2 u\circ\phi_{\xi,\lambda} \ \ \mathrm{for}\ \lambda\in(0,1] \ \ \mathrm{and} \ \ F(\xi,0):=\lim_{\lambda\downarrow 0^+}F(\xi,\lambda).\\[5pt]
\end{equation*}
The map $F$ is obviously continuous with $F(\xi,0)=u(\xi)$ for every $\xi\in\s$, i.e. $F(\s,0)=u(\s)=\s$, whereas $F(\s,1)=\{b_u\}$.
In other words, $F$ is a continuous homotopy between $\s$ and the point $b_u\in \overline {B_1}\diagdown \{0\}$ and therefore there exists $\lambda_0\in(0,1)$ such that $0\in F(\s,\lambda_0)$, i.e. there exists $\xi_0\in \s$ such that $F(\xi_0,\lambda_0)=0$.\\[-12pt]

In the general case of a map $u\in \mathcal{A}_{\s}$, by the approximation property given in Lemma 7, Section I.4. in \cite{brezis1995degree} there exists a sequence $(u_j)_{j\in \mathbb{N}}\in C^{\infty}(\s;\s)$ with the property that
\begin{equation*}\label{approximation_property}
u_j\underset{j\to \infty}{\longrightarrow} u \mathrm{\ strongly \ in} \ W^{1,2}(\s;\s) \mathrm{\ \ and \ \ } \mathrm{deg}u_j=\mathrm{deg}u=1\ \ \forall j\in \mathbb{N}.
\end{equation*} 
Up to passing to a non-relabeled subsequence, we can without loss of generality also suppose that $u_j\rightarrow u$ and $\t u_j\rightarrow \t u$ pointwise $\mathcal{H}^{2}$-a.e. on $\s$, as $j\to\infty$. Since the maps $u_j$ are smooth and surjective, by the previous argument there exist $(\xi_j)_{j\in \mathbb{N}}\in \s$ and $(\lambda_j)_{j\in \mathbb{N}}\in (0,1]$ so that for every $j\in \mathbb{N}$,
\begin{equation*}\label{sequence_of_Möbius_fixing_centers}
\d2 u_j\circ \phi_{\xi_j,\lambda_j}=0.
\end{equation*}

Up to non-relabeled subsequences we can suppose further that $\xi_j\rightarrow\xi_0\in \s$ and $\lambda_j\rightarrow\lambda_0\in [0,1]$ as $j\to \infty$, thus $\phi_{\xi_j,\lambda_j}\to \phi_{\xi_0,\lambda_0}$ pointwise $\mathcal{H}^{2}$-a.e. on $\s$ and also weakly in $W^{1,2}(\s;\s)$.\\[-12pt]

In fact $\lambda_0\in (0,1]$, i.e. the Möbius transformations $(\phi_{\xi_j,\lambda_j})_{j\in \mathbb{N}}$ do not converge to the trivial map $\phi_{\xi_0,0}(x)\equiv\xi_0$. Indeed, suppose that this was the case. Then $u_j\circ\phi_{\xi_j,\lambda_j}\to u(\xi_0)$ pointwise $\mathcal{H}^{2}$-a.e. and $|u_j\circ\phi_{\xi_j,\lambda_j}|\equiv1$, so we could use the Dominated Convergence Theorem to infer that
\begin{align*}
\begin{split}
u(\xi_0)=\d2 u(\xi_0)=\lim_{j\to\infty}\d2 u_j\circ\phi_{\xi_j,\lambda_j}=0, \\[5pt]
|u(\xi_0)|=\d2 |u(\xi_0)|=\lim_{j\to\infty}\d2 |u_j\circ\phi_{\xi_j,\lambda_j}|=1,
\end{split}
\end{align*}
and derive a contradiction. Having justfied that $\lambda_0\in (0,1] $, what we actually obtain by the Dominated Convergence Theorem is that 
\begin{equation}
\d2 u\circ \phi_{\xi_0,\lambda_0}= 0. 
\end{equation}
Hence, if we set $\psi:= \phi_{\xi_0,\lambda_0}$ and $\tilde u:=u\circ\psi$, thanks to the invariance of the Dirichlet energy and the degree under orientation-preserving conformal reparametrizations of $\s$, we have
\begin{equation}\label{properties_of_tilde_u}
\tilde u\in \mathcal{A}_{\s} \mathrm{\ \ with \ \ } D(\tilde u)=D(u), \ \ \mathrm{deg}\tilde u=\mathrm{deg}u=1\ \ \mathrm{ and\ \ } \d2 \tilde u=0.
\end{equation}
The proof is now divided in two further steps, where \textit{we first prove a local version of the Theorem under the assumption that our map $\tilde u$ is apriori sufficiently close to the $\mathrm{id}_{\s}$ in the $W^{1,2}$-topology, and then we use a contradiction$\backslash$compactness argument to reduce to the previous situation.}\\[-10pt]
	
{\bf Step 2}. Let us first  prove \eqref{conf_s2_s2_quantitative} under the extra assumption that 
\begin{equation}\label{conformal_s2_close_to_identity}
\d2 |\t \tilde u-P_T|^2\leq \theta^2,
\end{equation}	
where $P_T:=\t \mathrm{id}_{\s}$ and $0<\theta\ll 1$ will be suitably chosen later. This assumption implies a trivial upper bound for the conformal deficit, since
\begin{equation}\label{trivial_upper_bound_for_the_conformal_deficit}
D(u)=D(\tilde u)\leq\d2 \big(|\t \tilde u-P_T|^2+|P_T|^2\big)-1\leq 1+\theta^2.\\[5pt] 
\end{equation}
Note that since $\d2 \tilde u=0$ and $|\tilde u|\equiv1$ $\mathcal{H}^2-$a.e. on $\s$, we have
\begin{align}\label{conformal_deficit_is_Poincare_deficit}
\begin{split}
D(\tilde u)&=\frac{1}{2}\d2|\t \tilde u|^2-\d2 |\tilde u|^2=\frac{1}{2}\d2|\t \tilde u-\nabla \tilde u_h(0)P_T|^2-\d2 |\tilde{u}-\nabla \tilde u_h(0)x|^2\\[6pt]
&\geq \left(\frac{1}{2}-\frac{1}{6}\right)\d2|\t \tilde u-\nabla\tilde u_h(0)P_T|^2,
\end{split}
\end{align}
i.e.	
\begin{equation}\label{first_stability_estimate_conformal_s2}
\d2 |\t \tilde u-\nabla \tilde u_h(0)P_T|^2\leq 3D(\tilde u)=3D(u).\\[3pt]
\end{equation}
Here $\tilde u_h:\overline{B_1}\mapsto \R^3$ is the harmonic continuation of $\tilde u$ in the interior of the unit ball, taken componentwise.
One obtains \eqref{first_stability_estimate_conformal_s2} by viewing extrinsically $\tilde u$ as a map from $\s$ to $\R^3$ and expanding it in spherical harmonics. To be more specific, for every $k\in\mathbb{N}$ let $H_{k}$ be the subspace of $L^2(\s;\R^3)$ consisting of vector fields whose components are all $k$-th order spherical harmonics, i.e. eigenfunctions of $-\Delta_{\s}$ corresponding to the eigenvalue $\lambda_k:=k(k+1)$. As it is well known, $L^2(\s;\R^3)=\bigoplus_{k=0}^\infty H_{k}$, the orthogonal sum being taken with respect to the $L^2$-inner product. If $\Pi_{k}$ denotes the corresponding orthogonal projection on $H_k$, we have that $\Pi_{0}\tilde u=\d2 \tilde u=0$ and as a general fact $\Pi_{1} \tilde u=\nabla \tilde u_h(0)x$. Since the first non-trivial eigenvalue of $-\Delta_{\s}$ is $\lambda_{1}=2$, with the coordinate functions $(x_i)_{i=1,2,3}$ spanning its eigenspace, the first line in \eqref{conformal_deficit_is_Poincare_deficit} follows from the $L^2$-orthogonal decomposition $u=\Pi_1u+(u-\Pi_1u)$. The subsequent inequality follows by the sharp Poincare inequality on $\s$ for maps with vanishing linear part, since the optimal constant in the latter is the second nontrivial eigenvalue of $-\Delta_{\s}$, i.e. $\lambda_{2}=6$.\\[5pt] 
\textit{Therefore, we only have to justify why in \eqref{first_stability_estimate_conformal_s2} one can replace $A:=\nabla\tilde u_h(0)$ with a matrix $R\in SO(3)$}, up to changing the value of the constant on the right hand side.\\[-12pt] 

Towards this end, observe that by the mean-value property and a basic $L^2$-estimate for harmonic functions in $\overline{B}_1$ (which could be proven for example by expanding $\tilde u$ in spherical harmonics and $\tilde u_h$ in the corresponding homogeneous harmonic polynomials),  \eqref{conformal_s2_close_to_identity} implies that
\begin{equation*}\label{nabla_at_0_almost_id}
|A-I_3|^2=\left|\dashint_{B_1} \nabla \tilde u_h-I_3\right|^2\leq \dashint_{B_1} |\nabla \tilde u_h-I_3|^2\leq\frac{3}{2}\d2|\t\tilde u-P_T|^2\leq \frac{3\theta^2}{2},
\end{equation*}
and by choosing $\theta>0$ sufficiently small, we can take $A$ to be invertible and such that 
\begin{equation}\label{matrixinversedeterminantestimates}
|A|^2, |A^{-1}|^2\ \in [2,4]\ \ \ \mathrm{and \ \ \ } \mathrm{det}A\in \left[\frac{1}{2},\frac{3}{2}\right].
\end{equation}
	
By the polar decomposition, $A=R_0\sqrt{A^tA}$ with $R_0\in SO(3)$ and if we label the eigenvalues of $\sqrt{A^tA}$ as $0<\alpha_1\leq \alpha_2\leq \alpha_3$ and set $\lambda_i:=\alpha_i-1$, $\lambda:=\lambda_1+\lambda_2+\lambda_3$ and $\Lambda^2:=\lambda_1^2+\lambda_2^2+\lambda_3^2$, we have
	
\begin{equation}\label{L2_small_conformal_s_2}
\Lambda^2=\mathrm{dist}^2\big(A;SO(3)\big)\leq|A-I_3|^2\leq\frac{3\theta^2}{2}.
\end{equation}
If we also set 
\begin{equation}
\tilde w(x):=A^{-1}\big(\tilde u(x)-Ax\big),
\end{equation}
we can use the fact that $\mathrm{deg}(\tilde u)=1$ and expand the integral formula for the degree around the $\mathrm{id}_{\s}$ to obtain
\begin{equation}\label{first_equation_by_using_the_degree}
\mathrm{det}A\left(1+Q_{V_3}(\tilde w)+\d2\big\langle \tilde w,\partial_{\tau_1}\tilde w\wedge\partial_{\tau_2}\tilde w\big\rangle\right)=1,
\end{equation}
where (referring for more details on the calculations to \cite{zemas2020}, Proof of Lemma 4.5, in particular formulas (4.9) and (4.10) therein)
\begin{equation}\label{quadratic_form_for_volume}
Q_{V_3}(\tilde w):=\frac{3}{2}\ds \Big\langle \tilde w, (\mathrm{div}_{\s}\tilde w)x-\sum_{j=1}^3x_j\t\tilde w^j\Big\rangle. 
\end{equation}
On the other hand, by expressing $\mathrm{det}A$ as a polynomial in the eigenvalues we have 
\begin{equation*}\label{det_tilde_u_betas}
\mathrm{det}A=1+\lambda+\frac{1}{2}(\lambda^2-\Lambda^2)+\lambda_1\lambda_2\lambda_3
\end{equation*}
and therefore the identity \eqref{first_equation_by_using_the_degree} can be rewritten as
\begin{equation}\label{intermediate_identity}
\frac{\Lambda^2}{2}=\left(\lambda+\frac{\lambda^2}{2}\right)+\lambda_1\lambda_2\lambda_3+\mathrm{det}A\left(Q_{V_3}(\tilde w)+\d2\big\langle \tilde w,\partial_{\tau_1}\tilde w\wedge\partial_{\tau_2}\tilde w\big\rangle
\right).
\end{equation}
This last identity actually leads to the desired estimate, i.e. from it we can conclude that
\begin{equation}\label{linear_part_close_to_SO(n)_conformal_s2_s2}
\mathrm{dist}^2\big(A;SO(3)\big)=\Lambda^2\leq c_1D(u)
\end{equation}
for a universal constant $c_1>0$.\\[5pt] 
Indeed, to estimate the summand $\left(\lambda+\frac{\lambda^2}{2}\right)$ notice that
\begin{align*}
|A|^2\leq \frac{3}{2}\d2|\t \tilde u|^2=3+3D(u)\implies \lambda\leq \frac{3}{2}D(u)-\frac{\Lambda^2}{2}\leq\frac{3}{2}D(u).
\end{align*}
We can now distinguish between the case $\lambda\geq 0$, where 

\begin{equation}\label{first_summand_case_1}
\lambda+\frac{\lambda^2}{2}\leq \frac{3}{2}D(u)+\frac{9}{8}[D(u)]^2 \overset{\eqref{trivial_upper_bound_for_the_conformal_deficit}}{\leq} \left(\frac{3}{2}+\frac{9}{8}\big(1+\theta^2\big)\right)D(u),
\end{equation} 
and the case $\lambda< 0$, where due to \eqref{L2_small_conformal_s_2}
\begin{equation}\label{second_summand_case_2}
\lambda+\frac{\lambda^2}{2}\leq \left(1-\frac{3}{2\sqrt2}\theta\right)\lambda <0.
\end{equation} 
For the second summand in the right hand side of \eqref{intermediate_identity}, by the arithmetic mean-geometric mean inequality and \eqref{L2_small_conformal_s_2} we can easily estimate	
\begin{equation}\label{second_summand}
\lambda_1\lambda_2\lambda_3\leq\left(\frac{\Lambda^2}{3}\right)^{\frac{3}{2}}\leq \frac{\theta}{3\sqrt 2} \Lambda^2.
\end{equation}
Regarding the last summands, the quadratic term in the expansion of the degree can be easily estimated by using the Cauchy-Schwarz and the Poincare inequality as
\begin{align}\label{Qvol_bound}
\begin{split}
|Q_{V_3}(\tilde w)|&\leq \frac{3}{2}\left(\d2|\tilde w|^2\right)^{\frac{1}{2}}\left(\d2 |(\mathrm{div}_{\s}\tilde w)x|^2+\Big|\sum_{j=1}^3x_j\t \tilde w^j\Big|^2 \right)^{\frac{1}{2}}\\[1pt]
&\leq \frac{3}{2}\left(\frac{1}{6}\d2|\t \tilde w|^2\right)^{\frac{1}{2}}\left(\d2 |\t \tilde w:P_T|^2+\Big(\sum_{j=1}^3x_j^2\Big)\Big(\sum_{j=1}^3|\t \tilde w^j|^2 \Big)\right)^{\frac{1}{2}}\\[5pt]
&\leq \frac{3}{2\sqrt{2}}\d2 |\t \tilde w|^2\leq\frac{3|A^{-1}|^2}{2\sqrt{2}}\d2 \big|\t \tilde u-AP_T\big|^2\leq \frac{18}{\sqrt 2}D(u),
\end{split}
\end{align}
the last inequality following from \eqref{first_stability_estimate_conformal_s2} and \eqref{matrixinversedeterminantestimates}.\\[5pt]
For the last term, by \textit{Wente's isoperimetric inequality} for $\tilde w$,   \eqref{trivial_upper_bound_for_the_conformal_deficit}, \eqref{first_stability_estimate_conformal_s2} and \eqref{matrixinversedeterminantestimates}, we can estimate
\begin{align}\label{degtildewestimate}
\begin{split}
\left|\d2\langle \tilde w,\partial_{\tau_1}\tilde w\wedge\partial_{\tau_2}\tilde w\rangle \right|& \leq \left(\frac{1}{2}\d2|\t \tilde w|^2\right)^{\frac{3}{2}}\leq \frac{|A^{-1}|^3}{2\sqrt 2}\left(\d2 \big|\t \tilde u-AP_T\big|^2\right)^{\frac{3}{2}} \\[7pt]
& \leq  \frac{4\sqrt{27}}{\sqrt 2}\big[D(u)\big]^{\frac{3}{2}}\leq\frac{4\sqrt{27(1+\theta^2)}}{\sqrt 2}D(u).
\end{split}
\end{align}
By plugging \eqref{first_summand_case_1}-\eqref{degtildewestimate} into the identity \eqref{intermediate_identity}, keeping in mind \eqref{matrixinversedeterminantestimates} and rearranging terms, we obtain
\begin{align*}\label{almostfinalestimate}
\left(\frac{1}{2}-\frac{\theta}{3\sqrt 2}\right)\Lambda^2\leq\frac{3}{2}\left(1+\frac{3}{4}(1+\theta^2)+\frac{18+4\sqrt{27(1+\theta^2)}}{\sqrt{2}}\right)D(u),
\end{align*}
which is precisely \eqref{linear_part_close_to_SO(n)_conformal_s2_s2}, after choosing $\theta>0$ sufficiently small to additionally satisfy for instance
\begin{equation*}
\frac{1}{2}-\frac{\theta}{3\sqrt 2}\geq\frac{1}{4}\ \ \ \mathrm{and \ \ then\ \ } c_1:=6\left(1+\frac{3}{4}(1+\theta^2)+\frac{18+4\sqrt{27(1+\theta^2)}}{\sqrt{2}}\right).
\end{equation*} 
Therefore, by combining \eqref{first_stability_estimate_conformal_s2} and \eqref{linear_part_close_to_SO(n)_conformal_s2_s2} and by the conformal invariance of the Dirichlet energy on $\s$,
\begin{equation}
\d2 |\t \tilde u-R_0P_T|^2\leq 2\d2 \left|\t \tilde u-\nabla \tilde u_h(0)P_T\right|^2+\frac{2}{3}\Lambda^2  \implies \d2 |\t u-\t \phi|^2\leq cD(u),
\end{equation}
where $\phi:=R_0\psi^{-1}\in Conf_+(\s)$ and $c:=6+\frac{2}{3}c_1>0$.\\[-3pt]
	
{\bf Step 3}. Arguing by contradiction, suppose that the statement of  \hyperref[main_thm]{Theorem 1.1.} is false. Then for every $k\in \mathbb{N}$ there exists a map $u_k\in \mathcal{A}_{\s}$ with $D(u_k)>0$ such that 
\begin{equation}\label{contradictory_estimate}
\d2 |\t u_k-\t\phi|^2\geq kD(u_k) \ \mathrm{\ for \ all\ }\phi\in Conf_+(\s).
\end{equation}
In particular, for $\phi\in Conf_+(\s)$ which we can fix for the following computation and for $k\geq 5$,
\begin{align*}
\begin{split}
&kD(u_k)\leq \d2 |\t u_k-\t\phi|^2=\d2 \left|\t (u_k\circ\phi^{-1})-P_T\right|^2\\[5pt]
&\ \ \ \ \ \ \ \ \ \leq 2 \d2 \left(|\t (u_k\circ\phi^{-1})|^2+|P_T|^2\right)=2\d2|\t u_k|^2+4=4D(u_k)+8\\[5pt]
\implies&D(u_k)\leq \frac{8}{k-4}. 
\end{split}
\end{align*}
By letting $k\to\infty$ we obtain $\lim_{k\to\infty}D(u_k)=0$. We can then use the compactness result  of Lemma A.2. in the end of Appendix A of \cite{zemas2020} for the case $n=3$ (see also similar results in the references therein) to obtain a contradiction.\\[5pt]
Indeed, what this compactness result guarantees, is that up to passing to a subsequence we can find $\psi_k\in Conf_+(\s)$ and $R\in SO(3)$ so that the maps $v_k:=u_k\circ\psi_k\in \mathcal{A}_{\s}$ satisfy
\begin{equation*}
\d2 v_k=0\ \ \mathrm{and \ \ \ } v_k\to R\mathrm{id}_{\s} \ \ \mathrm{strongly\ in\ }  W^{1,2}(\s;\s) \ \ \mathrm{as\ \ } j\to \infty.
\end{equation*} 
Without loss of generality (up to considering $R^tv_k$ instead of $v_k$ if necessary) we can also suppose that $R=I_3$. Then, for the dimensional constant $\theta$ chosen in Step 2, we can find $k_0:=k_0(\theta)\in \mathbb{N}$ such that 
\begin{equation*}
\d2 |\t v_k-P_T|^2\leq \theta^2 \ \ \forall  k\geq k_0.
\end{equation*}
In other words, after precomposing with the correct Möbius transformations and also rotating properly, the subsequence  $(v_k)_{k\geq k_0}$ satisfies the condition $\d2 v_k=0$ and also fulfills the apriori closeness to the identity assumption \eqref{conformal_s2_close_to_identity}. By Step 2 we deduce that there exist $(\phi_k)_{k\geq k_0}\in Conf_+(\s)$ such that
\begin{equation*}\label{final_contradiction}
\ds |\t v_k-\t \phi_k|^2\leq c D(v_k) \ \ \forall k\geq k_0.
\end{equation*}
Combining now the last estimate with \eqref{contradictory_estimate} we arrive at the desired contradiction.  \qedhere

\begin{remark}
Interestingly, this first observation that the Möbius group of $\s$ can be used to fix the mean value of maps in $\mathcal{A}_{\s}$ to $0$ was also used in the proof of \hyperref[main_thm]{Theorem 1.1.} by P. Topping in \cite{topping2020}. While here we used it essentially to link the problem with the stability of the sharp Poincare inequality on $\s$, Topping uses it to start then the harmonic map heat flow with initial datum a map $u\in \mathcal{A}_{\s}$ with $\d2 u=0$. With this centering, the flow does not produce a bubble in finite time and converges as $t\to\infty$ to an orientation-preserving Möbius map of $\s$. This limiting map turns then out to be the one for which the desired stability estimate is satisfied. It would be interesting to see if these two approaches could be compared and further linked. 
\end{remark}

\section*{Acknowledgements}
J.H. is supported by the German Research Foundation (DFG) in context
of the Priority Program SPP 2026 ``Geometry at Infinity''. K.Z. is also supported by DFG  under Germany’s Excellence Strategy EXC 2044 -390685587, ``Mathematics Münster: Dynamics-Geometry-Structure''.` The content of this note was also included in K.Z.'s PhD Thesis, which was carried out at the Max Planck Institute for Mathematics in the Sciences in Leipzig, and was submitted to the University of Leipzig in June 2020.

\vspace{1em}


\vspace{3em}
$^1$ Universität Leipzig, Mathematisches Institut, Augustusplatz 10, 04109 Leipzig, Germany\\ 
\hspace{1em}\textit{Email address:} hirsch.jonas@math.uni-leipzig.de
\\[-3pt] 

$^2$ Applied Mathematics Münster, University of Münster, Einsteinstrasse 62, 48149 Münster, Germany\\ 
\hspace{1em} \textit{Email address:} konstantinos.zemas@uni-muenster.de 
\end{document}